\documentclass[11pt]{amsart}

\usepackage{lipsum}
\usepackage{mwe}
\usepackage{amsfonts}
\usepackage{amssymb}
\usepackage{mathrsfs}
\usepackage{latexsym}
\usepackage{graphicx}
\usepackage{fancyhdr}
\usepackage[english]{babel}
\usepackage{amsmath}
\usepackage{amsthm}
\usepackage{enumerate}
\usepackage{color}
\usepackage{cases}
\usepackage[lofdepth,lotdepth,caption=false]{subfig}
\usepackage{enumitem}
\usepackage{xcolor}
\usepackage{mathptmx}
\usepackage{float}
\usepackage{dsfont}
\usepackage{csquotes}

\usepackage{hyperref}


\usepackage[ruled,vlined]{algorithm2e}

\newcommand{\R}{\mathbb{R}}

\makeatletter
\def\blfootnote{\xdef\@thefnmark{}\@footnotetext}
\makeatother

\begin{document}
\title[AA-PINNs]{Auto-Adaptive PINNs with Applications to Phase Transitions}


 \author[Buck]{Kevin Buck$^{1,*}$}
 \author[Kim]{Woojeong Kim$^1$}
 
 \blfootnote{$^1$Institute for Scientific Computing and Applied Mathematics, Indiana University, Bloomington 47405 \\
 $^*$Corresponding Author.
 {\textit kevbuck@iu.edu} (K.~Buck), {\textit wki1@iu.edu} (W.~Kim)
 }

\begin{abstract}
We propose an adaptive sampling method for the training of Physics Informed Neural Networks (PINNs) which allows for sampling based on an arbitrary problem-specific heuristic which may depend on the network and its gradients. In particular we focus our analysis on the Allen-Cahn equations, attempting to accurately resolve the characteristic interfacial regions using a PINN without any post-hoc resampling. In experiments, we show the effectiveness of these methods over residual-adaptive frameworks.
\end{abstract}


\maketitle

\section{Introduction}
The study of Physics-Informed Neural Networks (PINNs) has grown rapidly in the past several years \cite{review,review2}. One key issue in the study of these objects is the difference between applications to stationary and time-dependent problems. Despite theoretical similarities, consistent methods for the training of PINNs on time dependent problems remain elusive. We propose that choosing efficient and intelligent sampling distributions are potentially a key to alleviating these issues.

In the simulation of statics with PINNs it is a common technique to oversample, undersample, or otherwise highlight in regions which are known to be problematic such as regions containing shocks \cite{Shock1,Shock2}, difficult boundary conditions \cite{BoundaryPINN}, or other irregularities \cite{Javi}. This is similar in spirit to adaptive mesh refinement schemes used extensively in traditional finite element methods, reviewed in \cite{AMR1,AMR2}. In this context it is well understood that choosing a proper mesh can make the difference between a problem being completely infeasible and easily solvable. Generally in static PINNs, problem regions are identified post-hoc and then manually sampled more heavily. This allows for the gradual improvement of the static function over the network training. As the training continues it is paused periodically, at which point the sampling distributions can be edited as needed. We refer to this methodology as post-hoc sampling, as it is done manually after some amount of training or detailed analysis of the particular approximation is done.

This post-hoc sampling is built on the fact that it is common in complex problems that some regions in the domain are more difficult to simulate than others. Simpler regions may have well-bounded error even without exceptionally low residual value, while other regions can explode in error with even the slightest inaccuracy. This is known as the conditioning of a problem \cite{NLA,me}. If small residual error leads to small actual accuracy, the problem is called well-conditioned. If a small residual does not necessarily result in small accuracy, the problem is called ill or poorly conditioned. More precisely, one can often bound the actual accuracy by some power of the loss times a constant. This constant is called the condition number. If it is very large, then the problem is poorly conditioned. 

In time-dependent problems, there are several additional difficulties when compared to statics. First is that the conditioning of a problem can vary in time and space. It could be that only one small region of a domain is poorly conditioned while the rest is generally well behaved. In this case the entire problem would be classified as poorly conditioned. Additionally the problem regions can move more drastically in the training process than in statics. This is due to the ability of the regions to move in time as well as in space combined with the local nature of the training of a PINN. During the training process, as earlier times shift to fit the PDE more thoroughly, the training already done on later times can be rendered irrelevant. This can also be true in space, as drastic shifting of one region of a spatial domain may affect the global solution of the PDE. These are key issues to address as they greatly diminish the effectiveness of training, especially if the network is somewhat far from the true solution.

To alleviate these issues, several methods have been experimented with. The first of these is Extended PINNs or XPINNs, first presented in \cite{xPINN}. This method involves the decomposition of the domain into non-overlapping subdomains. A separate PINN is then trained on each of these subdomains with an additional interior boundary condition guaranteeing that the PINNs match within the domain. Additional work has been done with this method to extend to variational problems \cite{varxPINN} and to determine which problems the method is effective in solving \cite{whenxPINN}. In a similar spirit \cite{WZ} proposes the partitioning of the time variable in their time slicing II approach, which will be a focus of our study.

This family of methods has several issues that we hope to address. First, the regions must be fixed at the start of training, and cannot be adjusted afterwards without restarting the entire training process. This is problematic if you expect the problem regions to move as the training is iterated. Additionally, XPINNs and its variants requires the training of multiple networks. The training of a PINN is extremely computationally expensive, and while the training of multiple PINNs can be done in parallel it may still be undesirable. Essentially this method is effective at splitting the difficulty in the domain into subdomains but does not directly address the issue of particularly problematic regions forming.

These ideas have motivated the broad studying of adaptive sampling methods, where training is enhanced through carefully examining the choice of sample points in network training. The first of these methods is known as residual adaptive (or loss adaptive) sampling, first introduced in \cite{orig_resadapt}, this method is explored thoroughly in \cite{resadapt}. This involves the training of a PINN by constantly resampling according to the pointwise residual of the network. The precise nature of how this sampling is done varies, but the key idea is that the loss is reduced more effectively by highlighting regions of greater loss by sampling them more heavily. This method is extremely widely applied, but has the theoretical shortcoming that it assumes a uniform loss is optimal. For problems which are highly spatially irregular, or have regions where error grows at vastly different rates, these problems become more obvious. Of particular interest here is that the method is easily applied without prior knowledge of where the problematic regions will be, a strength that we hope to maintain in our own methodology.

Additional strategies in adaptive sampling are discussed in \cite{R3,adapt1,adapt2,adapt3}. In particular, failure-informed PINNs (FI-PINNs) \cite{FIPINN,FIPINN2} present a flexible adaptive sampling framework which could in principle be combined with a methodology similar to our Energy-Adaptive PINNs. However, sampling strategies based on distributions which are not residual-based are in general underexplored in the literature.


We propose a training method aimed at alleviating the described issues by adaptively changing the sampling distribution used in training based not on the residual but on other network-dependent heuristics which indicate regions of high error and high error growth. For our study we focus on adaptively sampling in space, though this methodology could additionally be extended to allow resampling in time. In particular we propose adaptive sampling according to arbitrary densities dependent on the network, i.e. the approximated PDE solution, and its derivatives. To do this we employ a Metropolis-Hastings sampling in parallel with the network training. This functionally allows for an ideal sampling distribution to be learned by training. Importantly this type of sampling can also be done without manual intervention during the training process, instead automatically adjusting the sampling distribution according to the pre-programmed heuristic. This allows for the sampling distribution to gradually follow problematic regions as they form, move, and change in the entire domain of the problem. We call this method Auto-Adaptive both because of its self-referential nature and self-sufficiency once training begins. We discuss the particulars of this method in Section \ref{s:AASampling}.

Our test case for these methods is the Allen-Cahn equation, which describes the phase separation of two fluids over time \cite{AC}. It does this through the evolution of the order parameter, which indicates the relative concentration of the two mixed fluids. We discuss these equations in more detail in section \ref{AllenCahnSection}. In particular we experiment with an energy adaptive variant of the described auto adaptive sampling, where we sample in proportion to the pointwise energy density of the Allen-Cahn equation. As we show heuristically in \ref{HeuristicSection}, regions of high energy often correspond well with the problematic regions described earlier. These regions contain behavior where a small residual can invite a disproportionate amount of error. Additionally these regions are often areas more important to the continued accuracy as time is evolved. 

This setting highlights the issues that we hope to address with our adaptive methods. Allen-Cahn displays a separation of time and spatial scales as the interfaces quickly form and slowly move. However, these interfaces only partially describe where problematic regions occur. In particular, we contrast our results with those of \cite{WZ}, who also experiment with Allen-Cahn and the related Cahn-Hilliard model. The authors pioneer the time-slicing method by testing on the Allen-Cahn and Cahn-Hilliard systems. In their experiments, they encounter several issues that we hope to address. First, they observe consistent large error on the interfaces formed in the Allen-Cahn evolution as well as other regions where the function is away from $\pm 1$. We will show that the energy adaptive approach greatly alleviates this issue. Additionally for problems which have particularly large separation of scales, their ``time-adaptive I" approach fails and they are forced to instead train several networks on the same problem in what they call the ``Time-Adaptive II" approach. We find this solution undesirable, as the training of multiple networks may be expensive.

In the remaining introduction, we will describe some background about Physics-Informed Neural Networks and some important properties of Gradient Flow Systems and in particular Allen-Cahn. In section 2 we present the Auto-Adaptive PINN. In section 3 we discuss the practical implementation in preparation for our numerical experiments. Section 4 contains the results from several numerical experiments comparing these methods to baseline PINNs and to the well-tested Loss-Adaptive methods. Finally, in Section 5 we draw conclusions and expand upon directions for future work.

\subsection{Physics Informed Neural Networks}
First introduced in \cite{PINN}, a Physics-Informed Neural Network is a method of approximating a partial differential equation (PDE) using a Neural Network structure with artificial `data' derived from the PDE. In this framework, we define a neural network which is itself an approximation to the solution to a PDE. We then iteratively update this approximation through gradient descent on a loss function, which measures the fidelity of the approximation to the governing equations. 

In particular, given a PDE of the form
\begin{gather} \label{Introduction:GenericPDE}
\begin{cases}{}
  LHS[u] = RHS(x) \text{ in }\Omega\times (0,T) \\
  u(0,x) = u_0(x) \text{ in }\Omega \\
  u(x,t) = g(x,t) \text{ on }\partial\Omega\times (0,T),
\end{cases}
\end{gather}
the corresponding loss function is given by the sum of the residual norms of each equation:
\begin{equation} \label{Introduction:eq:PINN_Loss}
  L(\theta) = ||LHS[u_\theta(x, t)]-RHS(x)||^2_{\Omega\times(0,T)} + ||u_\theta(0,x) - u_0(x)||^2_\Omega + ||u_\theta-g||^2_{\partial\Omega\times (0,T)} .
\end{equation}

Here $\Omega\subseteq \mathbb{R}^d$ and $LHS$ represents the left hand side of the PDE, a differential operator on $u$. Meanwhile $RHS(x, t)$ represents the right hand side of the PDE, a forcing term dependent only on $x$ and $t$. The second equation in \ref{Introduction:GenericPDE} represents the initial condition and the third equation represents the boundary condition. For simplicity we write a Dirichlet type boundary condition but Neumann, mixed, or any other boundary condition can be represented similarly. 

Meanwhile in the Loss Equation \eqref{Introduction:eq:PINN_Loss}, $u_\theta(x, t)$ is the network itself, which depends on hidden parameters $\theta$ in a composite nonlinear way and takes inputs in $\Omega\times[0,T)$. The norms represent appropriate norms on the function spaces, which can vary depending on the problem. Typically, these norms are chosen to be $L^2$ norms on their respective function spaces. This yields an approximate loss function which is the mean square error of the residuals on the domain. These norms are then approximated via Monte-Carlo methods during the gradient descent process, yielding an approximate (or empirical) loss function that is used for computation. This Monte-Carlo Sampling will become essential for our adaptive methods, as we will substitute it for other more complex sampling distributions. This will amount to switching the usual $L^2$ norm to a weighted $L^2$ norm with weights dependent on the function itself, which will be discussed in detail in section \ref{s:AASampling}.

\subsection{Allen Cahn and Gradient Flow Systems}
We are motivated in particular by the Allen-Cahn equation. We will discuss first Gradient Flow problems in generality and then the particular choice of energy which leads to Allen-Cahn.

\subsubsection{Gradient Flow Equations}
A Gradient Flow Equation is one of the form
\begin{equation}\label{Introduction:eq:Gradient_Flow}
  \partial_t u = -\nabla_u E[u]
\end{equation}
where $E: H\to\R$ is coercive on $H$ a Hilbert space. We assume additionally that the energy has local structure,
\begin{equation}\label{Introduction:eq:Energy_Form}
  E[u] = \int_\Omega e(u, \nabla u) dx
\end{equation}
\noindent for some functional $e$. We then recall that for gradient flow systems we have the following equality:
\begin{equation}\label{Introduction:eq:Energy_Evolution}
  \frac{d}{dt} E[u(t)] = \langle\nabla E[u(t)], \partial_t u\rangle_H =- ||\partial_t u(t)||_H^2,
\end{equation}
which tells us exactly the rate of energy decay.

\subsubsection{Allen-Cahn Equation} \label{AllenCahnSection}
For our numerical examples we use the Allen-Cahn equation, which represent the $L^2$ gradient flow equation associated with the Ginzberg-Landau Free Energy. We take the Ginzberg Landau Free Energy
\begin{equation} \label{Introduction:eq:GinzbergLandauEnergy}
  E[u] = \gamma_1\int_\Omega |\nabla u|^2 dx + \gamma_2\int_\Omega \Psi(u) dx
\end{equation}
where $\Psi(s)$ is the Free-Energy Density and $\gamma_1$ and $\gamma_2$ are constants determining the relative energy balance. This free energy density is generally taken to be of a `double-well' type, meaning that is is minimized at exactly two points with one local max between the minimizing points. Thus the energy has two components, one which pushes the function towards two particular values (generally $\pm 1$), and one which penalizes steep gradients.

This energy yields the following as the Allen-Cahn equation.

\begin{equation}\label{Introduction:eq:AC_General}
  \partial_t u = \gamma_1\Delta u - \gamma_2\Psi'(u)
\end{equation}

For our experiments, we use the Landau Approximation of the Free-Energy Density, given by
\begin{equation}\label{Introduction:eq:Free_Energy_Density}
  \Psi(s) = \frac{1}{4} (s^2-1)^2.
\end{equation}
This has been shown as an appropriate choice of $\Psi$, though sometimes truncations are used to avoid numerical difficulties \cite{GinzLandApprox}. This yields the following final equation, which amount to a full description of our problem.

\begin{equation}\label{Introduction:eq:AC_Specific}
  u_t = \gamma_1 \Delta u - \gamma_2 (u^3-u)
\end{equation}

This equation describes the phase separation of a mixed fluid by modeling the evolution of their relative densities $u$, often called the order parameter. This quantity varies from $-1$ to $1$, with $-1$ representing purely one fluid while $1$ represents purely the other fluid. Generally it is assumed that $\gamma_1 << \gamma_2$, leading to a large separation of scales in the two terms. This generates the characteristic behavior of the system, which is to quickly evolve into regions close to $1$ and $-1$ with steep interfaces of characteristic width proportional to $\sqrt{\gamma_1}$. These interfaces then evolve on a slower timescale. 

\subsubsection{Error Heuristic for Allen-Cahn} \label{HeuristicSection}
In the case of the Allen–Cahn equation, we develop a heuristic for identifying regions of high error growth in PINN simulations, derived from well-established finite-difference and finite-element error estimation principles. This heuristic will be important for this implementation of Auto-Adaptive Sampling. We will study the growth of the error of a trained neural network $u_\theta$. 

First we define $\epsilon =u-u_\theta$ the error. We denote the network residual $\delta u$, which is given by
\begin{equation}
   \delta u := \delta_tu_\theta - \gamma_1\Delta u_\theta + \gamma_2 \Psi'(u_\theta)
\end{equation}
Then given $u$ an exact solution to \eqref{Introduction:eq:AC_General} we compute the PDE governing error growth as
\begin{equation}
  \delta_t\epsilon = \gamma_1\Delta\epsilon - \gamma_2(\Psi'(u)-\Psi'(u_\theta)) + \delta u.
\end{equation}
Linearizing the nonlinear term we get
\begin{equation}\label{Energy:equation:linearizederror}
  \delta_t\epsilon = \gamma_1\Delta\epsilon-\gamma_2\Psi''(u)\epsilon + \delta u.
\end{equation}
This linearized equation governs the growth of the error of the simulation with respect to time. If we assume that $\delta u$ is small, i.e. the approximation is close to accurate, we can see that this is a dispersive PDE with damping or amplification depending on the sign of $\Psi''(u)$.

Since we assume $\gamma_1 << \gamma_2$ and $\delta u$ is small, the $\Psi$ term \eqref{Introduction:eq:Free_Energy_Density} dominates the error grown equation \eqref{Energy:equation:linearizederror}. A low energy region is thus one where the value stays close to $\pm 1$ while a high energy region is one with values that stay close to $0$. If $u\approx \pm1$ then $\Psi''(u)\approx2$, and thus the corresponding term in \eqref{Energy:equation:linearizederror} becomes a dampening factor. Meanwhile in the high energy regions where $u\approx0$ we see that $\Psi''(u)\approx -1$, and the term corresponding to $\Psi''$ in \eqref{Energy:equation:linearizederror} becomes an amplification factor.

This heuristic indicates that the problematic regions in the simulation of Allen-Cahn will be areas of high energy. Thus we can see that regions of low energy naturally dampen the growth of the error while regions of high energy naturally amplify the growth of the error when the residual is reasonably small relative to the true solution.

\section{Auto-Adaptive Sampling}\label{s:AASampling}
In order to address the issues of moving problematic regions in time dependent PINNs, we introduce the Auto-Adaptive Sampling method. The premise of this method is to use a heuristic dependent on the network approximation itself, which indicates regions of high error growth, as a sampling distribution for the network training process. By sampling in proportion to this heuristic, we can directly reduce errors in the regions which are most susceptible and reach an optimal distribution of the residual loss.

Stated precisely, we replace the analytic formulation of the first term of \eqref{Introduction:eq:PINN_Loss} with a term of the following form
\begin{equation}\label{Energy_Adaptive_Sampling:eq:Weighted_Loss}
  L_{PDE}(\theta) = \int_0^T\int_{\Omega}|LHS[u_\theta(x,t)]-RHS(x)|^2 \rho(u_\theta) dxdt.
\end{equation}
Here $\rho(x)$ is taken proportionally to a constant $C$ plus a known heuristic function with $\int_\Omega\rho dx = 1$ and $\rho>0$ so that $\rho$ may be interpreted as a probability density (and any scaling can be seen as a reweighting of the PDE loss term). The addition of $C$ is important, as the function must be strictly positive to be a sensible sampling weight (there is no region which can be neglected entirely in training).

Rather than computing this weight directly, we will use it as a sampling distribution for the training process. By sampling in this manner and weighting all points equally in the loss computation, we implicitly recover the weight without explicitly computing its value. In addition, this naturally refines our distribution to be finer in difficult regions. In classical finite element schemes, this can be seen as analogous to a adaptive mesh refinement methods, which refine their mesh in the regions which are measured to be problematic.

This comparison is also why we decide to sample in proportion to the heuristic instead of use it as a directly computed weight on sample points. In low regularity phenomenon, the including of many sample points in the low regularity region is important to accurately capture the behavior. Using a smaller number of points with a directly computed weight does increase the emphasis of that region in training, but may fail to accurately capture the interior behavior of the region. The use of a sampling distribution is also beneficial for computational efficiency, which will be discussed in more detail with the Metropolis-Hastings Algorithm in section \ref{ss:MetropolisHastingsAlgorithm}.

With this framing, the constant $C$ can be reinterpreted as denoting the proportion of uniformly sampled points as opposed to points sampled in proportion to the heuristic. Thus we denote by $\lambda$ the ratio of points sampled adaptively:
\begin{equation}\label{Energy_Adapaptive_Sampling:eq:lambda_def}
  \lambda = \frac{\text{number of adaptive points}}{\text{number of adaptive points + number of uniform points}}
\end{equation}
then we compute the total PDE loss according to this ratio
\begin{equation}
  Loss_{PDE} = \lambda L_{adaptive} + (1-\lambda)L_{uniform}.
\end{equation}
The precise value of this hyperparameter will be chosen by experimentation. This decomposition also allows us to separately weight the importance of the adaptive and uniform points, which we will do in our numerical examples.

\subsection{Energy-Adaptive Sampling}\label{ss:EnergySampling}

For our simulations of the Allen-Cahn equation, we will use the network dependent heuristic provided in section \ref{HeuristicSection} as our sampling density. That is, we will sample in proportion to the pointwise energy density of the function approximation. 

Explicitly, our target density for the Metropolis Hastings Algorithm will be given (up to scaling) by:
\begin{equation}
  \rho_A(x) = \gamma_1 |\nabla u_\theta|^2 + \gamma_2\Psi(u_\theta), \hspace{.5cm} \rho_A(\cdot) \text{ is uniform in }t.
\end{equation}
where $\Psi$ is the double well potential defined by \eqref{Introduction:eq:Free_Energy_Density}. 

We observe that this density would be expensive to invert, thus motivating our use of the Metropolis-Hastings algorithm. With this algorithm, it is only necessary to compute the values of this density at our current and proposed sample points. Since gradients of the current sample points are already available to the network from the optimization algorithm, and the density values of the proposed sample points can be easily calculated via automatic differentiation, the algorithm can be efficiently employed in parallel to network training.

Notably, we impose that the density must be uniform time. This is because if allowed to move in time, sample points will coalesce at the beginning of the time interval (as the total energy of the system decreases according to equation \eqref{Introduction:eq:Energy_Evolution}). Exploring time dependent adaptive densities in an obvious direction to be explored in future work, and could be used to address complex time-dependent irregularities.

We then combine this with uniformly sampled points to get an effective total sampling density
\begin{equation}\label{eq:targetdisttribution}
  \rho(x,t) = C + \alpha \rho_A(x)
\end{equation}
where $C$ is defined implicitly through $\lambda$, the proportion of points uniformly sampled as defined by \eqref{Energy_Adapaptive_Sampling:eq:lambda_def}. We denote by $\alpha$ a hyperparameter accounting for the relative importance of the adaptively sampled points to the training. This can be interpreted as weighting the relative importance for the separation of scales. With this observation, we choose this weight as exactly $\alpha =\sqrt{\gamma_2/\gamma_1}$, which dictates the separation of scales.

Points in this sampling are then weighted equally in the discrete loss computation as to truly recover the weighted loss term \eqref{Energy_Adaptive_Sampling:eq:Weighted_Loss}.

\subsection{Metropolis-Hastings Algorithm} \label{ss:MetropolisHastingsAlgorithm}
In order to efficiently sample from a probability distribution dependent on the network and its derivatives, we use the Metropolis-Hastings Algorithm \cite{Metropolis,Hastings}. This is a sampling method which iteratively improves a collection of randomly chosen points to match the distribution given by a known probability density function. This method avoids the need to invert a distribution, and is thus ideal for settings where direct sampling is difficult. 

As a brief description, given some random sample of points $x$, the Metropolis-Hastings algorithm proposes new points, $x'$ drawn from some proposed probability density $g(x'|x)$ which can be easily sampled from and may depend on $x$. The points are then individually accepted or rejected probabilistically according to their adherence to the target probability density $\pi(x)$, along with a corrective ratio called the Metropolis Ratio. In particular a proposal point is accepted with probability
\begin{equation}\label{MH:correctiveRatio}
  A = \min\left(1, \frac{f(x')g(x'|x)}{f(x)g(x|x')} \right)
\end{equation}
where $f(x)$ is any function proportional to the target density function $\pi(x)$. Iterating the process many times yields sample points $x_N$ with a density approaching the target density as $N$ becomes large. This method has been shown to be effective even in situations with fairly pathological target density functions, given proper choice of $g$.

It is worth noting that there is a brief history of the use of Metropolis-Hastings and the broader class of Markov-Chain-Monte-Carlo (MCMC) methods in PINNs. In particular, MCMC methods are suggested by \cite{resadapt} as a method of implementation for residual adaptive methods when the sampling distribution is complex. Though dropout methods remain more popular in residual adaptive methods, such as was employed in \cite{WZ}, due to their ease of implementation. Additionally the Bayesian-PINN (B-PINN) \cite{BPINN} uses MCMC-based inference to capture uncertainty in the network parameters, addressing noise through probabilistic modeling rather than deterministic regularization, though the way the methods are used in this context is very different from our own.

The Metropolis-Hastings is additionally very well-studied and lends itself very naturally to the context of Physics-Informed Neural Networks. As an iteration based method, it can be run easily in conjunction with training. Additional improvements such as multiple proposal \cite{MultiMCMC} can nicely address common issues in PDEs and in shocks in particular. Additionally, the method is nicely parallelizable and can be run very efficiently on GPUs \cite{ParMCMC}. These benefits make it an exceptionally good fit for sampling in the context of Neural Networks.

In our context we will use the Metropolis-Hastings method not to sample according to a likelihood or posterior distribution as is typical, but to sample according to our network-dependent heuristic function $\rho$ described above in the above Section \ref{s:AASampling} and in \ref{ss:EnergySampling}.

\subsection{Details of Metropolis Hastings}
In this section we describe in more detail the precise nature of the Metropolis-Hastings Procedure we use in our experiments. When initializing the adaptive points on each time slice, we perform $1,000$ iterations of the Metropolis-Hastings Algorithm described in \ref{ss:MetropolisHastingsAlgorithm}. We use a proposal distribution which is a normal distribution centered around each point with standard deviation $\sqrt{\gamma_1}$, which is proportional to the interfacial width. We additionally update the standard deviation with each iteration to target an acceptance ratio (the proportion of proposed points accepted) of between $0.2$ and $0.6$.

The target distribution is given by the pointwise energy \eqref{eq:targetdisttribution}. This is easily computed by the network since the relevant gradients of existing sample points are available from the backpropagation used to compute the residual. The gradients of the proposed points can be computed in the same manner. The corrective ratio is given from the explicit formula for the pointwise density function of the normal distribution centered on each point. The acceptance probability \eqref{MH:correctiveRatio} is then computed directly from these values.

We also `step' the adaptive points at the conclusion of each training epoch (not after each minibatch). To do this, we perform $200$ iterations of the same algorithm described above. This ensures that the sampling distribution continuously changes to reflect changes to the solution approximation (i.e. the network).

\section{Other Techniques}\label{s:OtherTechniques}
In order to facilitate training, we use alongside the Auto-Adaptive method several standard techniques. In this section we discuss the details for the implementations of each of these methods. Details which fluctuate between examples are not specified here, but can be found in the detailed discussion of each example. At the conclusion of this detail-oriented discussion we supply an algorithm schematic which presents an overview of the methodology.

\subsection{Time Slicing}
This method was introduced in \cite{WZ} and involves the gradual expansion of the time domain in discrete steps. The idea is to reduce the time complexity of the problem by first greatly restricting the time domain until it is palatable to the network. Once the problem is learned on that smaller interval, the interval can be increased. This allows for the network to only require learning a small time interval at any given time time, which is generally easier for training. This method was tested very thoroughly with a large degree of success. However, for higher-dimensional and less regular systems further issues arise.

One of the these issues is commonly refereed to as `Catastrophic Unlearning' in the wider Machine Learning community, though use of this term is not commonly used in the discussion of PINNs. This issue is common in traditional data driven contexts where networks are trained to do one task and then trained to do a second separate task. In this instance the network may unlearn the original task in favor of learning the second task. In the context of PINNs this issue manifests when a network which is well trained on some sub-interval of the entire time domain, is trained on a separate or expanded time interval. As a particular example, a PINN trained on the time interval $[0,0.1]$ may drastically lose accuracy when trained on $[0,0.2]$ since simulating from 0.1 to 0.2 can be considered a `new task.' Not only will the simulation not perform accurately on the new time window of 0.1 to 0.2, but it will also lose the accuracy it had on the original window of 0 to 0.1. Of course simulating accurately on $[0,0.1]$ is essential to simulating on $[0.1,0.2]$, so in the setting of PINNs the problem is even further accentuated.

The authors of \cite{WZ} were aware of this issue, so they presented an alternative method that they call `Time Slicing II.' This method involves the training of a separate network on each time slice rather than using the gradual expansion of training time on a single network. This does resolve the issue of catastrophic unlearning, since no network is trained to perform multiple tasks. However, it is somewhat undesirable as it requires many networks, a true time discretization, and can also dramatically increase cost in training time and number of network parameters. For these reasons we manually edit the learning rate as we evolve through time slices, which is discussed next.

\subsection{Learning Rate Schedule}
We additionally use a learning rate scheduler, which updates the learning rate as the training process progresses. We use a fairly simple scheduler which updates the learning rate only based off the current time slice being trained. This helps to combat the catastrophic unlearning described in the section above. We decrease the learning rate roughly linearly as the progress through the time-slices. More sophisticated methods of choosing the learning rates are constantly being studied \cite{lrannealing,softadapt}, however for simplicity and ease of implementation and comparison we choose only to use this simple scheme.

\subsection{Residual Adaptive Sampling}
A commonly employed method to increase training efficacy is Residual Adaptive sampling, first proposed in \cite{orig_resadapt}. This technique involves sampling the domain in proportion to the current pointwise value of the residual. This leads to decreasing the loss very efficiently, as has been thoroughly explored in \cite{resadapt} and seen in many other papers including \cite{WZ}. In our numerical experiments, we will use this method as a test to compare against energy adaptive sampling. 

As was discussed in the introduction, the effectiveness of this method is generally built on the idea that a uniform loss value is desirable. In actuality, uniform loss does not necessarily minimize error. This is especially true in problems where the accuracy on a small region has a disproportionate effect on the overall error, such as problems with sharp interfaces or other separation of spatial scales. 

\subsection{Minibatching}
Rather than use the entire set of collocation points for each iteration of gradient descent, we instead use minibatching in order to achieve faster convergence. Minibatching is the practice of subsampling a larger collection of collocation points for each iteration. In particular, for a total number of sample points $N$ we randomly select only $N_{mini}<<N$ for each iteration. On each further iteration the points previously sampled are omitted until the entire original $N$ sample points are chosen or until fewer than $N_{mini}$ points remain. One pass of this procedure is called an \textit{epoch}. So the number of iterations per epoch is $\lfloor N/N_{mini}\rfloor$. Note that if $N$ is not divisible by $N_{mini}$ some points will be excluded.

\subsection{Latin Hypercube Sampling}
For uniform sampling in our domain, we use Latin Hypercube sampling. This method reduces variance in sampling by partitioning the domain into a grid of $N^d$ hypercubes, where $N$ is the desired number of sample points and $d$ is the number of spatial dimensions. The hypercubes are then sampled randomly so that the selected hypercubes are orthogonal. This yields a selection of $N$ hypercubes in the grid, with no selected hypercube in the same row or column as any other. Once this orthogonal set of hypercubes is selected, one point within each is chosen randomly from a uniform distribution on that hypercube. These individual sample points combine to form a representative sample of $N$ points from the domain. Further explanation can be found \cite{hypercube}. 

This method guarantees that the uniform sample will be representative of the entire domain. This sampling also greatly reduces variance, and thus is less reliant on the convergence yielded by the law of large numbers. This especially helps to reduce inconsistencies of sampling too few points, and makes the training more robust.

\subsection{Initial Condition Weight}
We use the common technique of enforcing a large weight on the initial condition term of the loss. In all cases we use a weight of $1000$ on the initial condition. This is to ensure that the initial condition is met as precisely as possible, since it is a hard constraint of our problem. Additionally, if the initial condition is not properly learned, the entire simulation is immediately rendered inaccurate.

\subsection{Optimizer Choice}
The primary method used in our numerical experiments is Adaptive Moment Estimation (ADAM). This method has been shown many times to work effectively in the context of neural networks and PINNs, and is a very standard choice. It is also well known that the Broyden-Fletcher-Goldfarg-Shanno (BFGS) method and related variants have been shown to work well when used after the application of ADAM, near the end of the training time. L-BFGS has two additional problems which should be addressed for our use case. First, should it be used at the end of each time-slice or only after all slices have been trained. Secondly, since BFGS requires fixed sample points over many iterations, exactly how do we intertwine this method with the adaptive sampling method. 

Answering these questions amounts to numerical experiment. We found that the implementation of BFGS at the conclusion of each time slice hurt overall results. We speculate that the method did too good a job at `locking in' the behavior of the system early on, which prevented the network from adapting from the fast behavior in the first portion of the time interval to the slow behavior in the later portion of the interval. For this reason we only apply ADAM for our results.

\subsection{Algorithm Schematic}
Combining the implementation methodology described above, we present the following algorithm schematic which consolidates the material.
\vspace{3cm}
\begin{center}\begin{algorithm}[H]
\SetAlgoLined
\KwResult{Train physics-informed neural network (PINN) for PDE}
\textbf{TrainingLoop}$(time\_slices, epochs, \lambda_{IC}, \lambda_{PDE})$: \\
\Indp
  Define domain and collocation sample sizes \\
  Compute $\lambda_{pde}$ and $\lambda_{IC}$ ratios \\
  \For{each $final\_time$ in $time\_slices$}{
    Set optimizer learning rate \\
    Initialize adaptive collocation points \\ 
    \For{epoch = 1 to epochs}{
    Sample uniform collocation points \\
      \For{each minibatch}{
        Compute loss (IC, PDE, BC terms with weights) \\
        Backpropagate gradients \\
        Apply gradient clipping \\
        ADAM Optimizer step
      }
      Update adaptive points \\
    }
  }
\caption{PINN Training with Auto-Adaptive Sampling}
\Indm
\end{algorithm}\end{center}

\section{Numerical Results}
The following numerical benchmarks are taken from \cite{WZ}. In the examples they provide, there is a clear concentration of the error in the regions of high energy. As such, this is an ideal setting for the testing of the improvements given by our new method.

For each example, we compare three methods. First, we run a high fidelity finite difference to yield an `exact' solution. The second method is Residual Adaptive sampling. Here we use the Metropolis Hastings Algorithm to sample in proportion to the loss and combine these with points sampled uniformly as suggested in \cite{resadapt} for complex sampling distributions. Additionally this implementation uses all of the techniques described in the implementation section, including time slices, learning rate scheduling, minibatching, latin hypercube sampling, and initial condition weights. Finally we implement our energy adaptive method described in section 3.

All simulations for this project are done in python using the PyTorch package for Neural Networks, and run on an NVIDIA Quatro RTX A6000 GPU, provided generously for use by the Institute for Scientific Computing and Applied Mathematics at Indiana University. All scripts can be found online at \href{https://github.com/kevmbuck/Auto-Adaptive-PINNs}{https://github.com/kevmbuck/Auto-Adaptive-PINNs}.

\subsection*{Example 1}
As our first example we use the parameters
\begin{equation}
  \gamma_1 = 1e-4,\hspace{.2cm} \gamma_2 = 5, \hspace{.2cm} u_0(x) = x^2\cos(\pi x)
\end{equation}
on the domain [-1, 1] with periodic boundary conditions.

\begin{figure}[H]
  \centering
  \includegraphics[width=.45\textwidth]{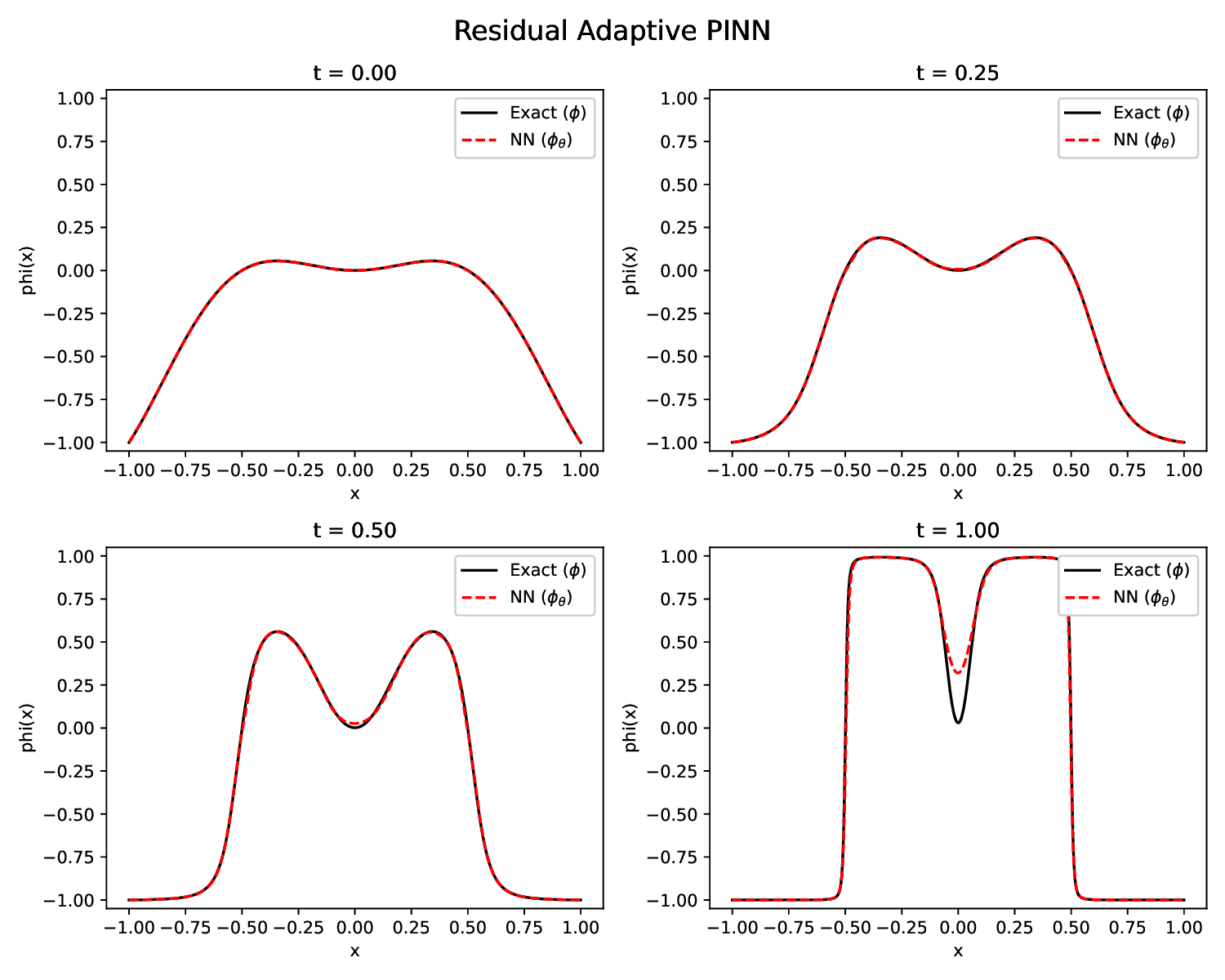}
  \includegraphics[width=.45\textwidth]{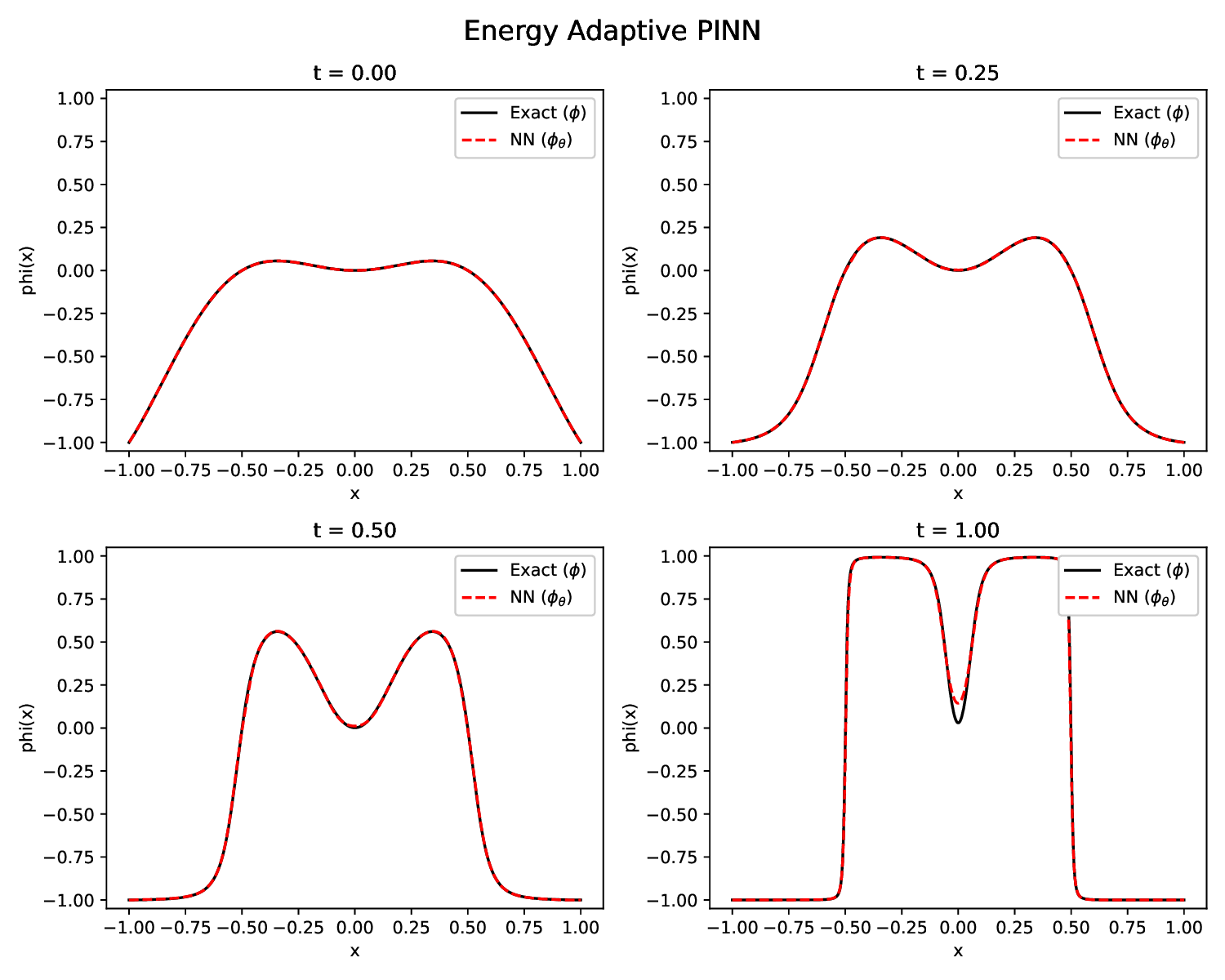}

  \caption{\small{Time slices from networks which have completed the described training. We present time slices from the the residual adaptive method and the energy adaptive method at times $0$, $0.25$, $0.5$, and $1$.}}
  \label{fig:Ex1:Slices}
\end{figure}

A region of interest in this example is the midpoint of the spatial domain, which can be seen in Figure \ref{fig:Ex1:Slices}. In this region although the function is far from one or minus one, the slow flow dominates since $u\equiv 0$ is an exact local maxima of the energy density. Additionally $\Delta u_0>0$ in the center, so the flow is nonzero. The dynamics at this point are thus very difficult, as any deviation from the exact state will change the flow dramatically.

Our network uses 6 fully connected layers of 128 nodes. We use the hyperbolic tangent activation function for all interior layers. For training we use $10,000$ total collocation points in the domain. We train the network until final time $t=1$, with time slices taken in increments of $0.1$. On each time slice, we run 100 epochs of ADAM training with minibatches of size $40$ split proportionally according to the adaptive point ratio, $\lambda$. As we progress through time slices, we additionally change the learning rate as described in section \ref{s:OtherTechniques}. We use learning rate $10^{-3}$ until final time 0.3, at which point we switch to $5*10^{-4}$. At final time $0.5$, we decrease to learning rate $10^{-4}$. At $0.7$ we decrease to $5*10^{-5}$, and finally at $0.9$ we decrease to $10^{-5}$.

\begin{figure}[H]
  \centering
  \includegraphics[width=0.975\linewidth]{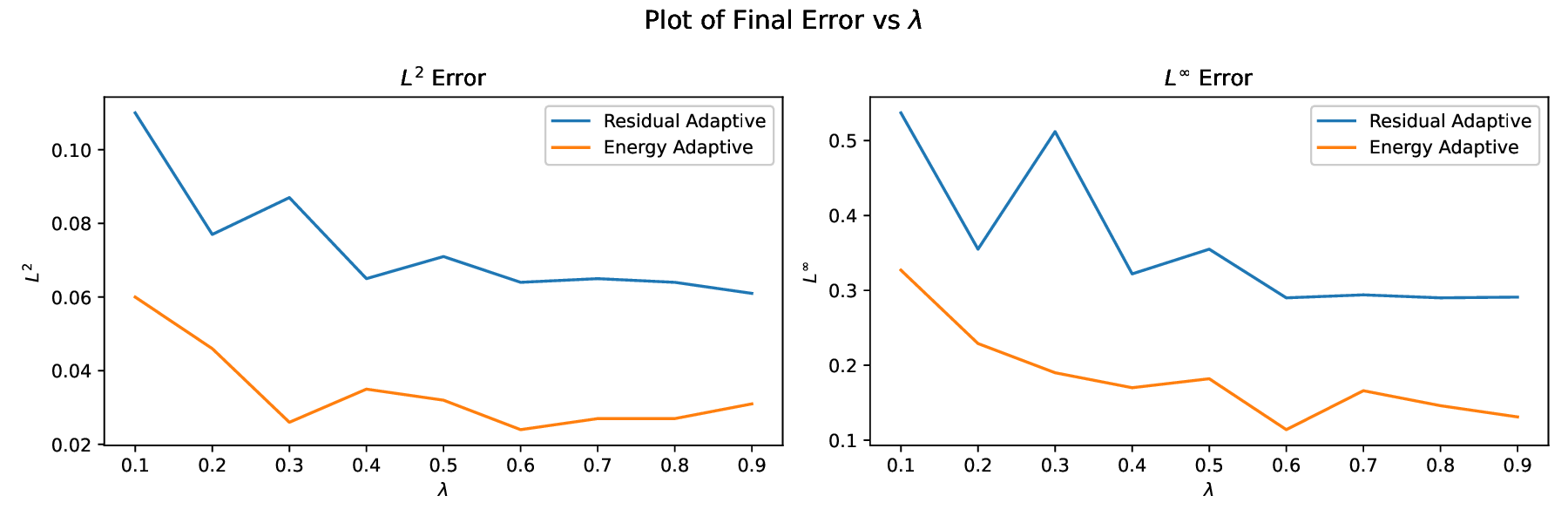}
  \caption{{\small Plots of the $L^2$ (left) and $L^\infty$ (right) error at the final simulation time as a function of $\lambda$, the proportion of adaptively sampled points. Each plot depicts the error of the residual adaptive in blue and the energy adaptive in orange. }}
  \label{fig:Ex1:lambda_sweep}
\end{figure}

We perform a parameter sweep on $\lambda$, the proportion of adaptively chosen points vs uniformly chosen points. The results for this are found in Figure \ref{fig:Ex1:lambda_sweep}. We can see from this plot that the energy adaptive method performs best with a $\lambda$ value of about $0.6$, while residual adaptive functions about equivalent for $\lambda$ between $0.6$ and $0.9$. For further discussion we consider $\lambda=0.6$ so as to more directly compare with the energy adaptive method. Notice also the consistent improvement of the energy adaptive method in both error measures.

\begin{figure}[H]
  \centering
  \includegraphics[width=.675\textwidth]{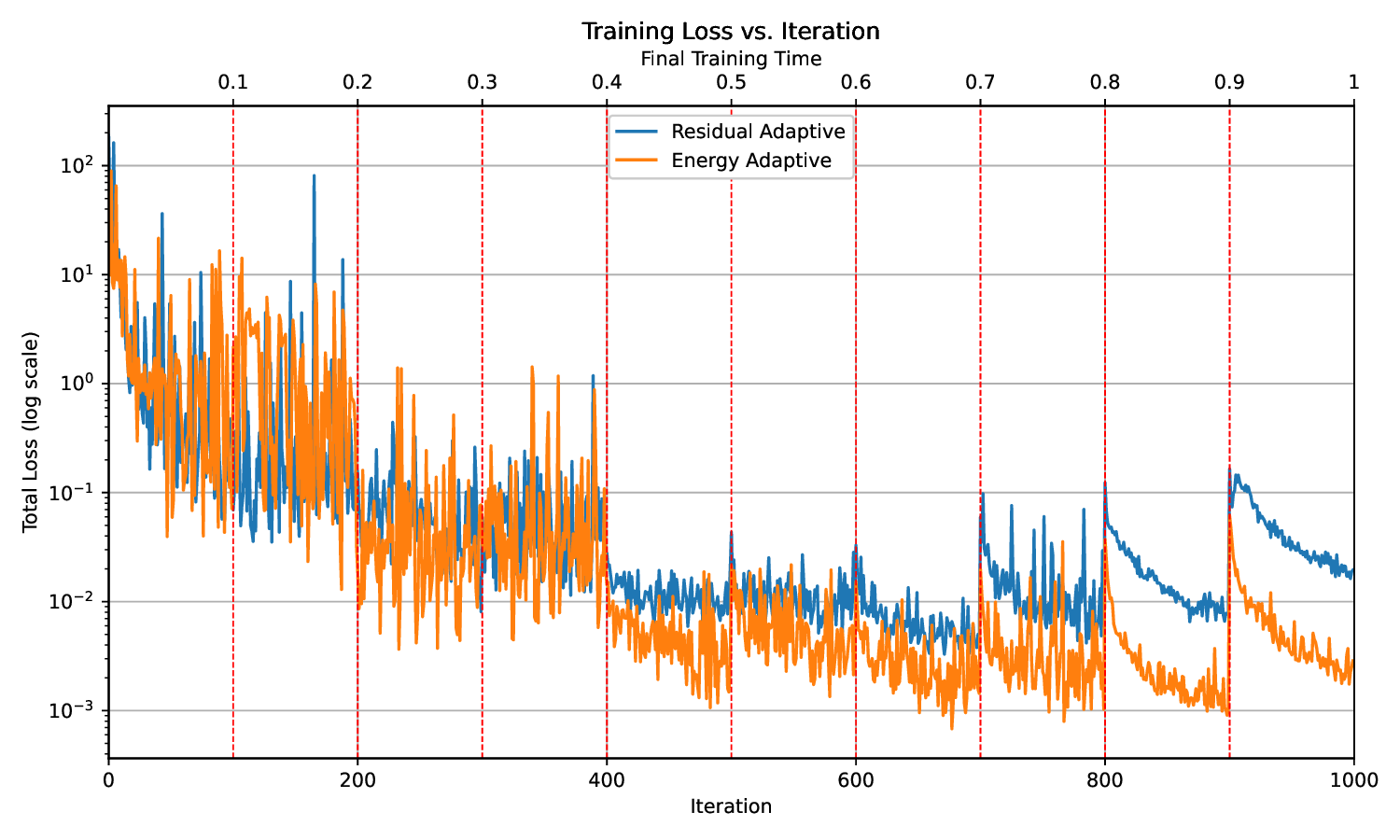}

  \caption{\small{Plot of the loss against the number of ADAM iterations performed. The vertical axis represents the base $10$ log of the loss value, while the horizontal axis represents the number of epochs trained. The vertical red lines represent the increasing of the trained time domain, which is labeled at the top of the graph. Note every other red line also corresponds with a reduction in learning rate.}}
  \label{fig:Ex1:Loss_History}
\end{figure}

To verify that both methods are decreasing loss effectively, we display the loss across each epoch in Figure \ref{fig:Ex1:Loss_History}. It is verified here that each method is successfully evolving according to their respective losses.

We next present time slices obtained from each of the tested methods in Figure \ref{fig:Ex1:Slices}. We notice the residual adaptive method generally performs well at early times but quickly loses accuracy in problematic regions as minor errors in the center expand dramatically. The energy adaptive method performs significantly better but it is not immune to this difficulty. We then observe the precise error measures in Table \ref{table:Ex1:Errors}, which verify exactly what we see in Figure \ref{fig:Ex1:Slices}.

\begin{table}[H]
  \centering 
  \begin{tabular}{| c | c | c | c |}  
  \hline
    Error Measure / Method & Residual Adaptive & Energy Adaptive \\ \hline
    Relative $L^2$ at $T=1$ & $4.09\times 10^{-2}$ & $1.50\times 10^{-2}$ \\ \hline
    $L^\infty(0, T; L^\infty(\Omega))$ Error & $2.31\times 10^{-1}$ & $7.96\times 10^{-2}$ \\\hline
  \end{tabular}
  \caption{\small{We observe the errors of each method. 
  }}
  \label{table:Ex1:Errors}
\end{table}

\begin{figure}[H]
  \centering
  \includegraphics[width=.975\textwidth]{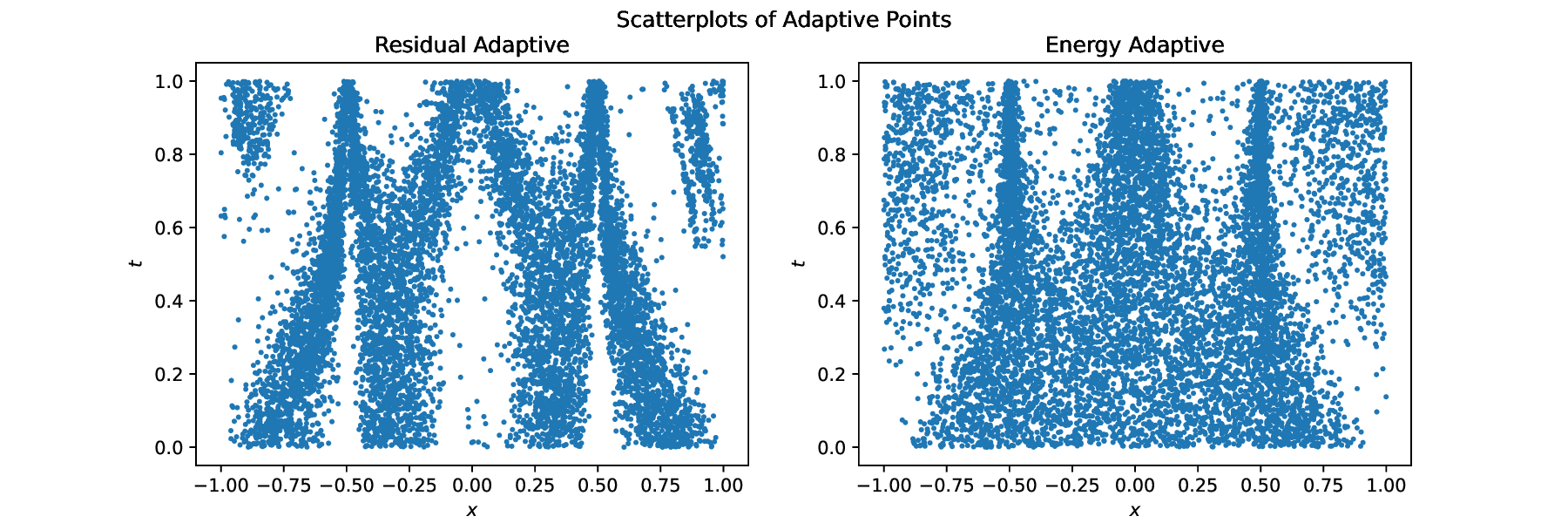}

  \caption{\small{Scatterplots containing the adaptively sampled points for the residual adaptive method (left) and the energy adaptive method (right). The horizontal axis represents the spatial domain and the vertical axis represents the temporal domain. These points are taken from near the end of network training and are thus representative of the entire spatiotemporal domain.}}
  \label{fig:Ex1:Scats}
\end{figure}

We additionally note here the success of the Metropolis Hastings method in capturing the distributions for both the residual and energy adaptive sampling in Figure \ref{fig:Ex1:Scats}. Notice in particular how the residual adaptive distribution captures areas surrounding the interfaces and center region while the energy adaptive method samples these regions directly. This captures the energy adaptive PINNs ability to \textit{directly} emphasize the problematic regions rather than \textit{indirectly} address them through the loss.

\subsubsection*{Computational Complexity}
We include some brief remarks on computational cost. Our residual-based and energy-based adaptive sampling methods each required approximately $5,500$ seconds of wall-clock time, which is expected since both rely on a Metropolis–Hastings procedure with comparable proposal and acceptance steps. By comparison, a baseline PINN trained with identical hyperparameters but without adaptive sampling required approximately $3,200$ seconds. Quick computations show that the baseline method with no adaptive sampling has a cost of $3.2$ seconds per epoch, while adaptive sampling incurs an additional cost of $2.3$ seconds per epoch and $4.6$ seconds per time-slice (the cost of re-initializing the adaptive points).

While adaptive sampling therefore incurs a nontrivial additional cost, the baseline PINN completely fails to capture the solution, whereas the adaptive methods succeed. Moreover, this cost compares favorably with dropout-based residual selection strategies, which typically introduce additional overhead due to repeated residual evaluation and resampling; in practice, the total training time often increases as the fraction of retained points decreases.

\subsection*{Example 2}
For our second example the following parameters determine the problem:
\begin{equation}
  \gamma_1 = 1\times 10^{-4},\hspace{.2cm} \gamma_2 = 4, \hspace{.2cm} u_0(x) = x^2\sin(2\pi x)
\end{equation}

This example is very similar to the above but we swap the even symmetry of the initial condition for an odd one. This makes the behavior more stable, as although the center region still has a local minima of the free energy density (the fast flow is $0$), since $\Delta u_0 = 0$ at the center the slow flow is also zero. Instead of monitoring the error in the upward drift as in Example 1, here we will monitor the ability of the methods to capture the split interface in the center of the domain. 

We additionally have problems with the boundary not seen in the first example, as there is an interface in the true solution going through the periodic boundary. Here the PDE loss at the boundary acts in competition to the boundary loss term. This yields to the PDE loss being much higher around the boundary as the two terms conflict with each other in the training process.

We use experimental hyperparameters all identical to the first example, in order to test the robustness of each method without fine tuning. The problems are also very similar, so these hyperparameters are not too far from optimal.

\begin{figure}[H]
  \centering
  \includegraphics[width=.7\textwidth]{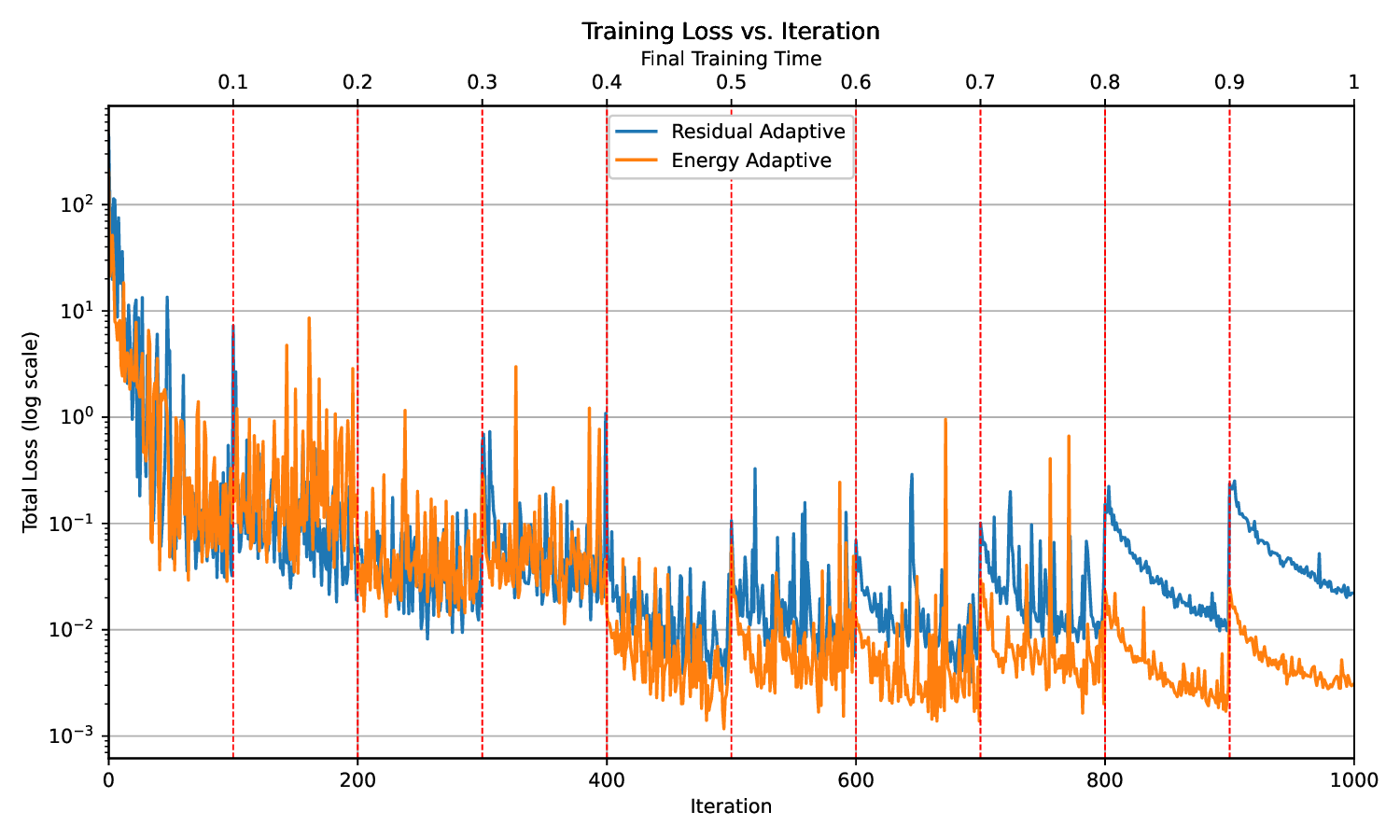}

  \caption{\small{We see the loss decreasing for each of the tested methods on the second example. The vertical axis represents the base $10$ log of the loss value, while the horizontal axis represents the number of epochs trained. The vertical red lines represent the increasing of the trained time domain, which is labeled at the top of the graph.}}
  \label{fig:Ex2:Loss_History}
\end{figure}

\begin{figure}[H]
  \centering
  \includegraphics[width=.475\textwidth]{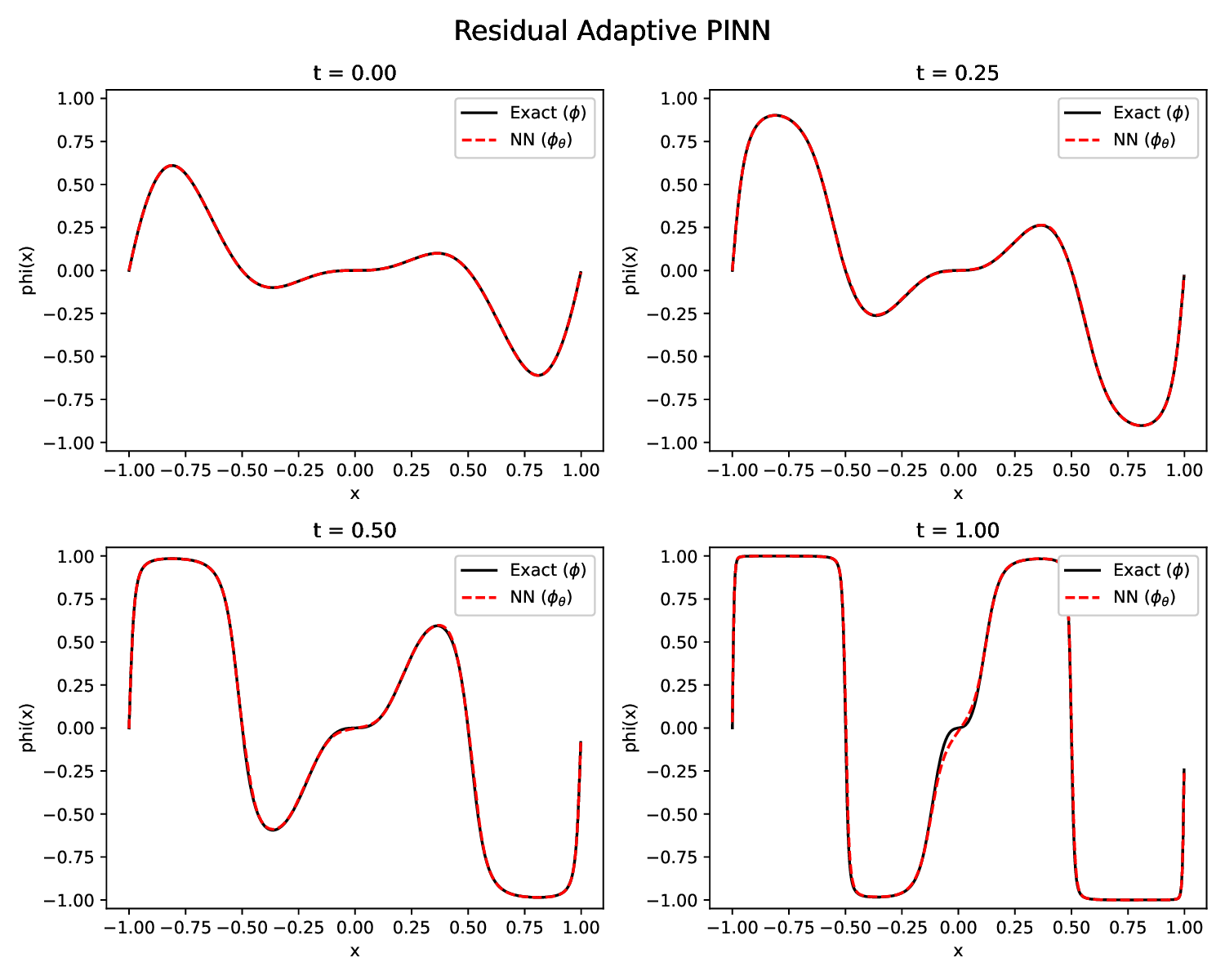}
  \includegraphics[width=.475\textwidth]{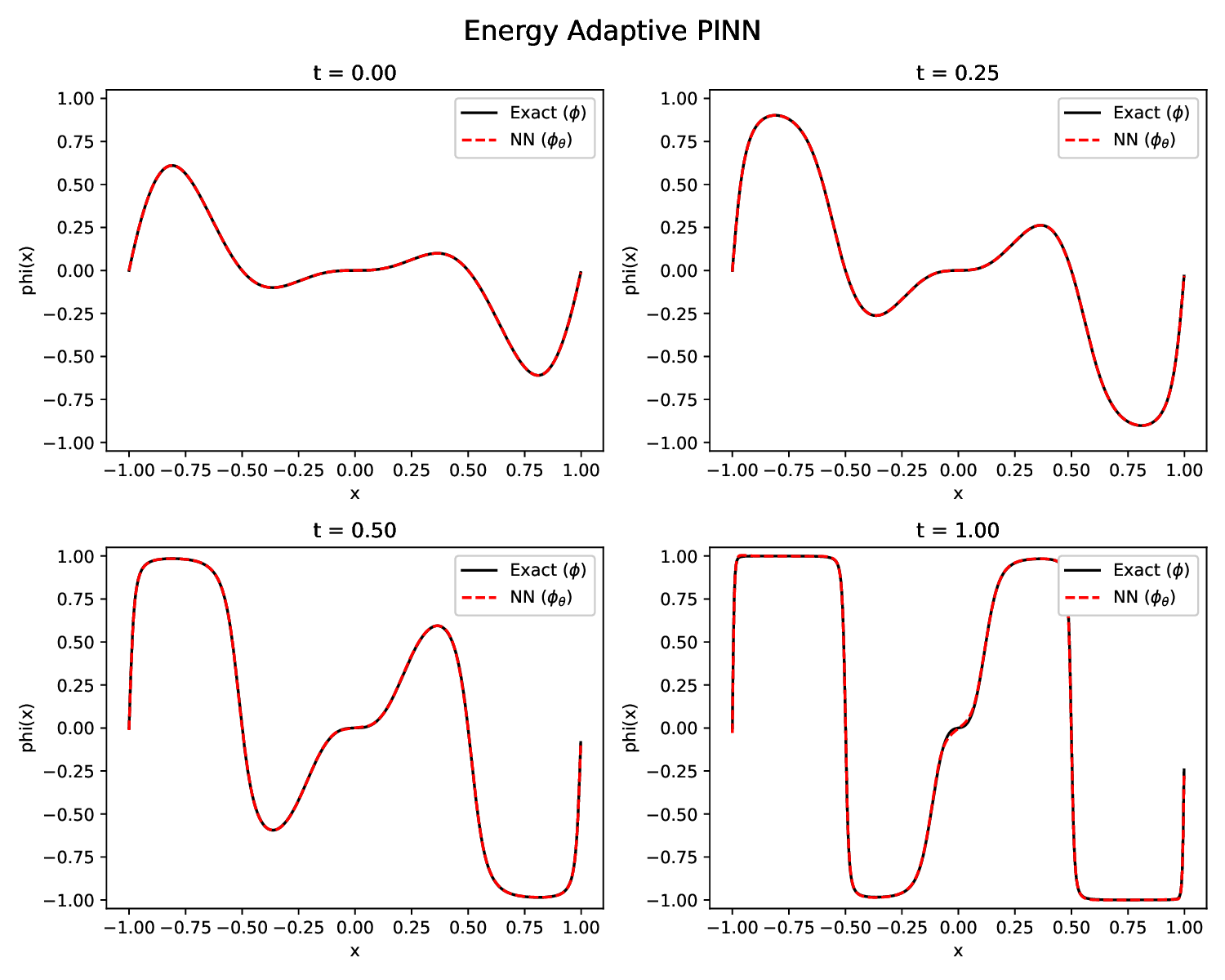}

  \caption{\small{Time slices from networks which have completed the described training. We present time slices from the residual adaptive method and the energy adaptive method at times $0$, $.25$, $.5$, and $1$. 
  }}
  \label{fig:Ex2:Slices}
\end{figure}

We again first observe that the loss is decreasing successfully for all methods throughout the training process in Figure \ref{fig:Ex2:Loss_History}. Then the results of these simulations are shown in Figure \ref{fig:Ex2:Slices}. We see again significant problems in the center region, though they do not result in error as high as the first example due to the lack of motion at the center. Instead, we can see the region losing the dual-interface structure of the exact solution in favor of a (residual-wise) simpler interpolant. Only in the energy adaptive method is the proper interfacial structure maintained. We additionally observe the relative $L^2$ and $L^\infty$ errors in Table \ref{table:Ex2:Errors} to verify the heuristics seen in the slices. We see here that the energy adaptive method improves over the Residual Adaptive method by an order of magnitude.

\begin{table}[H]
  \centering 
  \begin{tabular}{| c | c | c | c |}  
  \hline
    Error Measure / Method & Residual Adaptive & Energy Adaptive \\ \hline
    Relative $L^2$ at $T=1$ & $2.33\times 10^{-2}$ & $6.87\times10^{-3}$ \\ \hline
    $L^\infty(0, T; L^\infty(\Omega))$ Error & $1.15\times 10^{-1}$ & $3.20\times 10^{-2}$ \\\hline
  \end{tabular}
  \caption{\small{We observe the errors of each method. 
  }}
  \label{table:Ex2:Errors}
\end{table}

\subsection*{Example 3}
Finally we experiment with the 2D Example given in \cite{WZ}. This is characterized by the parameters
\begin{equation}
  \gamma_1 = \lambda\epsilon^2, \hspace{.3cm} \gamma_1 = \lambda, \hspace{.3cm} u_0(x) =\tanh\left(\frac{.35-\sqrt{(x-.5)^2+(y-.5)^2}}{2\epsilon}\right) 
\end{equation}
with $\lambda = 10$, $\epsilon = .025$, spatial domain $[0,1]\times[0,1]$, and final time $10$. This experiment is also different in that the initial condition is already very close to having the interfaces formed. Thus the slow behavior of the interfaces dominates the flow. Between times $.9$ and $1$ the behavior changes as the interfaces recede and the solution moves to the constant state at $u\equiv-1$. This late in time change of behavior presents an interesting environment to test our methods, in addition to the difficulties naturally presented by adding an additional spatial dimension.

\begin{figure}[H]
  \centering
  \includegraphics[width=.825\textwidth]{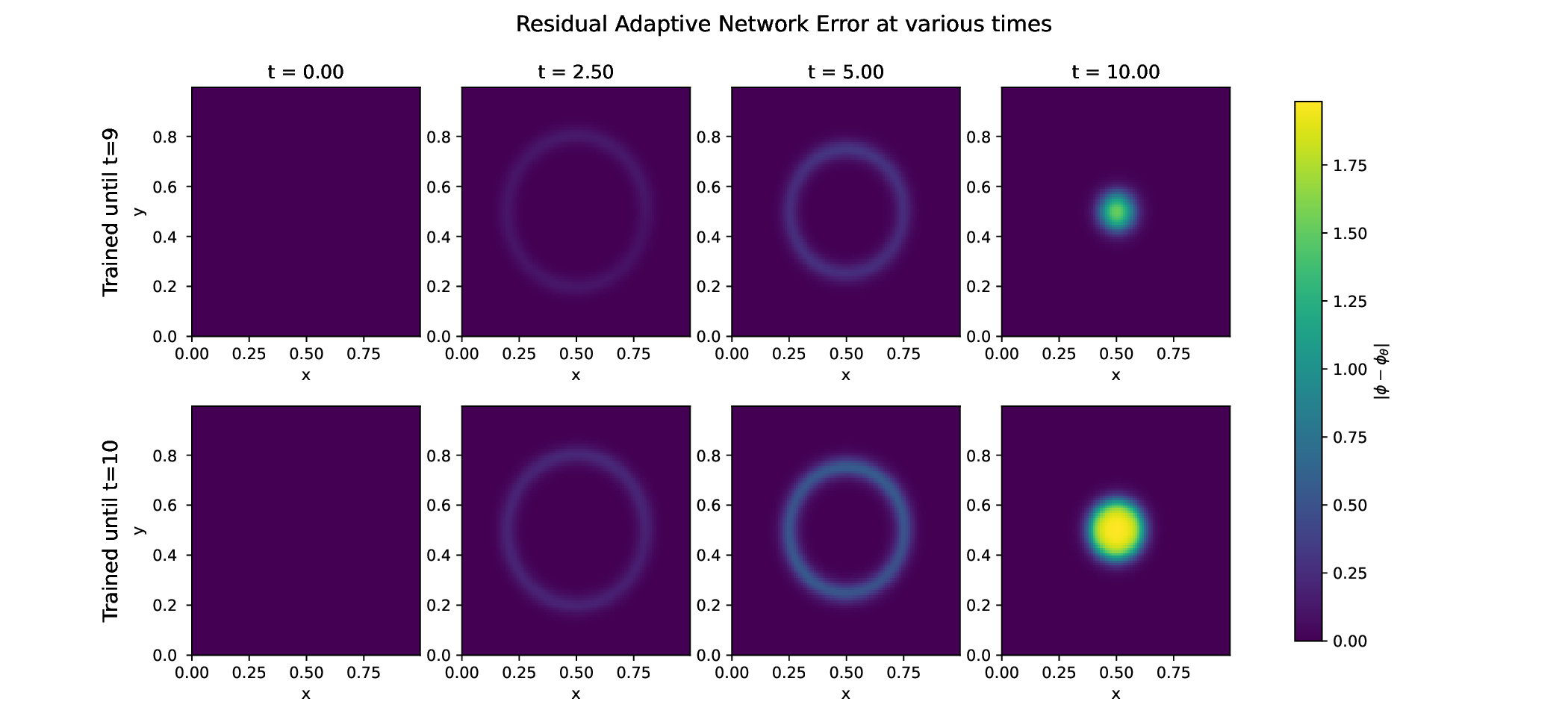}

  \caption{\small{We show the network error at various times. The top row shows the network only trained to a final time of $9$ (thus the depiction at $t= 10$ is the network attempting to extrapolate from what was already learned). The bottom row shows the results of the network after training until final time $10$.}}
  \label{fig:Ex3:Unlearning_residual}
\end{figure}

We can see this difficulty manifest in Figure \ref{fig:Ex3:Unlearning_residual}, which depicts results from the baseline residual adaptive method. Here we see that the expansion of the training time from $t=9$ to $t=10$ results in the catastrophic unlearning described in the introduction. Not only do we loose accuracy on the final time slice, but we also loose accuracy on the entire preceding time domain.

\begin{figure}[H]
  \centering
  \includegraphics[width=.825\textwidth]{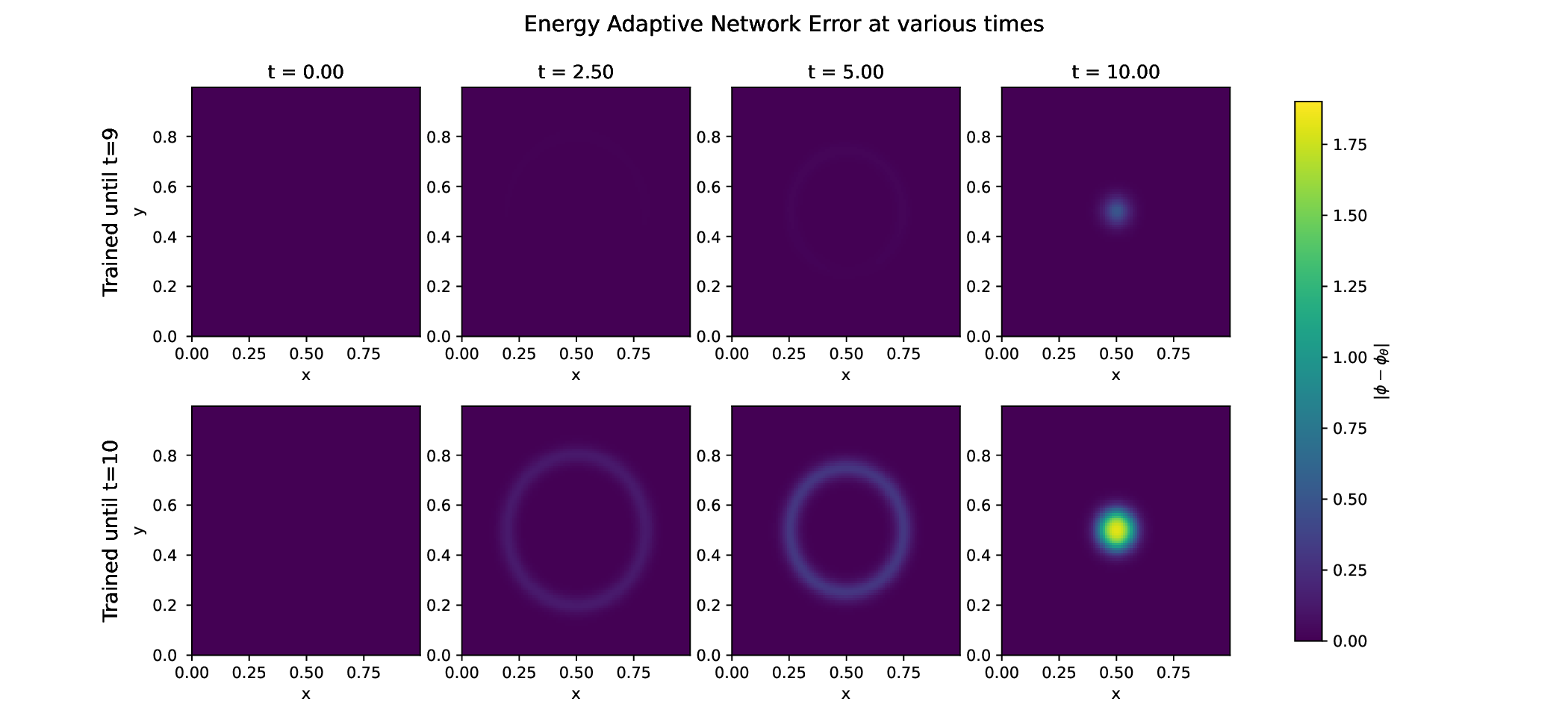}

  \caption{\small{We show the network error at various times. The top row shows the network only trained to a final time of $9$ (thus the depiction at $t= 10$ is the network attempting to extrapolate from what was already learned). The bottom row shows the results of the network after training until final time $10$.}}
  \label{fig:Ex3:Unlearning_energy}
\end{figure}

We next show the results from energy adaptive training in Figure \ref{fig:Ex3:Unlearning_energy}. We see here that the errors are much smaller and more contained than in the residual adaptive method. The precise error values can additionally be found in Table \ref{table:Ex3:Errors}. We see that the energy adaptive method when trained until time $9$ leads to very good extrapolation until time $10$. Here the errors are acceptable even for a simulation on the full domain. Despite the improved accuracy, the energy adaptive method still succumbs to the catastrophic unlearning phenomenon in the final time slice. Although even after the unlearning the errors are still comparable to the residual adaptive method before unlearning.

\begin{table}[H]
  \centering 
  \begin{tabular}{| c | c | c | c |}  
  \hline
    Error Measure / Method & Trained until & Residual Adaptive & Energy Adaptive \\ \hline
    Relative $L^2$ at $T=10$ & $T=9$ & $1.52\times 10^{-1}$ & $3.56\times 10^{-2}$ \\ \hline
     Relative $L^2$ at $T=10$ & $T=10$ & $3.53e-01$ & $1.93\times 10^{-1}$ \\ \hline
    $L^\infty(0, T; L^\infty(\Omega))$ Error & $T=9$ & $1.73\times 10^0$ & $9.79\times{-1}$ \\\hline
    $L^\infty(0, T; L^\infty(\Omega))$ Error & $T=10$ & $1.98\times 10^0$ & $1.90\times 10^0$ \\\hline
  \end{tabular}
  \caption{\small{We observe the errors of each method. 
  }}
  \label{table:Ex3:Errors}
\end{table}

In general there are several methods which can be used to reduce the problems posed by catastrophic unlearning. The two most prominent of which would be the employment of an xPINN \cite{xPINN} or the use of exponential time sampling \cite{ExpTime}. These directions are known to be effective, but require a post-hoc knowledge of where the unlearning occurs. For our analysis here, we are trying to confine our approach only to methods which can be applied without prior knowledge of the location of the unlearning. We also believe that an intelligently chosen heuristic employed via auto-adaptive PINNs could efficiently resolve this issue. However, we leave the exploration of this as a future direction.



\section{Conclusions and Future Directions}
In this study we have proposed an adaptive sampling method which allows for the training of Physics-Informed Neural Networks to be done with a complex sampling density dependent on the network and derivatives of the network itself. In particular we have shown the effectiveness of sampling in proportion to the energy of the system for the Allen-Cahn equations. This has allowed for the capture of complex dynamics with large separation of scales in both slow and fast timescales, and greatly alleviated the characteristic issue of these equations to concentrate a large amount of error into very small regions.


How the method behaves on a larger class of equations is still an open problem. In particular, more benchmarks for the performance of the energy adaptive method on Allen-Cahn and on the related Cahn-Hilliard system would be beneficial. Although we developed the energy adaptive method specifically for the Ginzberg-Landau energy in particular, we conjecture that its use on other gradient flow problems could alleviate issues for other complex choices of energy such as Total Variation Flow, Mean Curvature Flow, the Thin Film Equation, the Fokker-Planck equation, and Wasserstein Gradient flow.

Additionally, using the Metropolis Hastings in parallel with the network training should be experimented with for more densities in a much wider variety of contexts. This method allows for the complete customization of the sampling density for any individual problems. This framework is extremely flexible and could alleviate issues in many other difficult time dependent problems. Additionally, our current MCMC implementation is relatively inefficient as we use only the baseline method with no GPU implementation beyond the use of pytorch tensors. These methods have a wide variety of potential improvements that could drastically increase the time performance of our employed method.

Analytic verification in the form of a proof is also desirable for the energy adaptive method, perhaps in a similar manner to the proofs in \cite{ExpTime}. If well understood, this could also motivate choices of sampling densities for other difficult problems.

Finally, we suggest that some hybrid of the residual and energy adaptive methods could be desirable. The residual adaptive method is extremely efficient at decreasing the loss uniformly. The energy adaptive method succeeds by acknowledging that a uniform loss is not in itself desirable for minimal error. A hybrid approach could potentially outperform either method individually if well-executed.

\section*{Acknowledgments}
We would like to thank the Institute of Scientific Computing and Applied Mathematics for generously providing access to their computing resources, including the two NNVIDIA A6000 GPUs used for our computations. The work of author Kevin Buck was also partially funded by NSF Grants number DMS-2154387 and DMS-2206105.


\end{document}